\documentclass[preprint,10pt]{elsarticle}

\newtheorem{theorem}{Theorem}[section]
\newtheorem{lemma}{Lemma}[section]
\newtheorem{corollary}{Corollary}[section]
\newtheorem{proposition}{Proposition}[section]

\newtheorem{remark}{Remark}[section]
\newcommand{\ignore}[1]{}{}

\usepackage{amssymb}
\usepackage{amsmath}

\usepackage{color}
\usepackage{soul}
\usepackage[dvipsnames]{xcolor}

\usepackage[colorlinks=true, linkcolor=blue, citecolor=blue]{hyperref}

\def\1{{{\mbox{${\rm{1\negthinspace\negthinspace I}}$}}}}

\newcommand{\eref}[1]{(\ref{#1})}

\newcommand\beq{\begin{equation}}
\newcommand\eeq{\end{equation}}

\usepackage{ifthen}
\usepackage{xkeyval}
\usepackage{todonotes}
\begin{document}

\begin{frontmatter}

\title{Cram\'{e}r moderate deviations for a supercritical Galton-Watson process}
\author[cor1]{Paul Doukhan}
\author[cor2]{Xiequan Fan}
\author[cor3]{Zhi-Qiang  Gao}
\address[cor1]{CY University, AGM UMR 8088, site Saint-Martin,
	95000 Cergy-Pontoise, France}
\address[cor2]{Center for Applied Mathematics,
Tianjin University, Tianjin 300072,  China}
\address[cor3]{Laboratory of Mathematics and Complex Systems (Ministry of Education), School of Mathematical Sciences,\\ Beijing Normal University, Beijing 100875, China}


\begin{abstract}
Let $(Z_n)_{n\geq0}$ be a supercritical Galton-Watson process. The Lotka-Nagaev estimator $Z_{n+1}/Z_n$  is a common estimator for the offspring mean.
In this paper, we establish some Cram\'{e}r moderate deviation results for the Lotka-Nagaev estimator via a martingale method.
 Applications  to construction of confidence intervals are also given.
\end{abstract}

\begin{keyword} Cram\'{e}r moderate deviations;   Lotka-Nagaev estimator; offspring mean
\vspace{0.3cm}
\MSC primary   60G42; 60F10; 60E15;    secondary  62F03
\end{keyword}

\end{frontmatter}




\section{Introduction}
\setcounter{equation}{0}

Let $(X_i)_{i\geq1}$ be a sequence of independent and identically distributed (i.i.d.) random
variables with mean $0$ and positive variance $\sigma^2.$  Denote by $S_n =\sum_{i=1}^{n}X_i$
 the partial sums of $(X_i)_{i\geq1}$.  Assume $\mathbb{E}\exp\{ c_{0}|X_{1}|\}<\infty$ for some constant $c_{0}>0.$
Cram\'{e}r \cite{Cramer38} has established the following asymptotic moderate deviation expansion: for all $0 \leq x =  o(n^{1/2}), $
\begin{equation}
\Bigg|\ln \frac {\mathbb{P}(S_n> x\sigma\sqrt{n})} {1-\Phi(x)} \Bigg| =  O \bigg( \frac{1+x^3}{\sqrt{n}}\bigg) \ \ \mbox{as} \ \ n \rightarrow \infty,
\label{cramer001}
\end{equation}
where $\Phi(x)=\frac{1}{\sqrt{2\pi}}\int_{-\infty}^{x}\exp\{-t^2/2\}dt$ is the standard normal distribution.
The results of type (\ref{cramer001}) are usually called as Cram\'{e}r moderate deviations.
 After the seminal work of Cram\'{e}r,
a number of Cram\'{e}r moderate deviations have  been established for various settings. See, for instance, Linnik \cite{L61}  and \cite{F17}   for independent random variables, Fan, Grama and Liu \cite{FGL13}  for martingales (see also Puhalskii \cite{P97} for large deviation principles), Grama,  Liu and Miqueu  \cite{GLE17} and Fan,  Hu and  Liu \cite{FHL20} for a supercritical branching process in a random environment,  Beknazaryan,  Sang and Xiao \cite{BXY19} for random fields, and Fan \emph{et al.} \cite{FGLS19,FGLS20} for self-normalized type results.  In this paper,
we are going to establish Cram\'{e}r moderate deviations for a supercritical Galton-Watson process.

A Galton-Watson process is defined as follows
 \begin{equation} \label{GWP}
 Z_0=1,\ \ \ \ Z_{n+1}= \sum_{i=1}^{Z_n} X_{n,i}, \   \ \ \ \textrm{for } n \geq 0,
 \end{equation}
where $X_{n,i}$ is the offspring number of   the $i$-th individual of the generation $n$ and $Z_n$ stands for the total population of   the generation $n$.
Moreover,  $  (X_{n,i})_{i\geq 1} $ are independent of each other with a common  distribution law  $\mathbb{P}(X_{n,i} =k   )     = p_k,\   k \geq0,$
  and are also independent to $Z_n.$ Denote by $m$ the offspring mean of an individual, then  it holds  $$m=\mathbb{E}Z_1=\mathbb{E} X_{n,i} =\sum_{k=0}^\infty  k   p_k,\ \ \ n, i \geq 1.$$
 Denote by $v$ the standard variance of $Z_1$, then we have
\begin{eqnarray}\label{defv}
  \upsilon^2 =\mathbb{E} (Z_1-m)^2=\textrm{Var}( X_{n,i})=\textrm{Var}(Z_1).
\end{eqnarray}
To avoid triviality, we assume that  $v$ is positive.
The Lotka-Nagaev estimator $Z_{n+1}/Z_n$  is a common estimator for the offspring mean $m$.
 Throughout the paper, we  assume that $p_0=0$, then the Lotka-Nagaev estimator $Z_{n+1}/Z_n$ is well defined $\mathbb{P}$-a.s.
 Athreya \cite{A94} has established  large deviation rates
for the Lotka-Nagaev estimator with $Z_1$ satisfying Cram\'{e}r's condition. See also Ney and Vidyashankar \cite{NV03,NV04}
under the assumption that $\mathbb{P} (Z_1 \geq x) \sim a x^{1- \alpha}, x\rightarrow \infty,$ for two constants $\alpha>2$ and $a>0.$
Fleischmann and  Wachtel \cite{FW08} considered a generalization of the Lotka-Nagaev estimator $S_{Z_n}/Z_n,$
where $S_n=\sum_{i=1}^{n}X_i$ is  independent of $Z_n$ and $\mathbb{P} (Z_1 \geq x)=\mathbb{P} (X_1 \geq x) \sim a x^{ - \beta }, x\rightarrow \infty,$ for a constant  $\beta>2$. See also He \cite{H} when $X_1$ is in the domain of attraction of a stable law.
For the Galton-Watson processes with immigration, we refer to Liu and Zhang \cite{LZ06} and Li and Li \cite{LL21} for the rates of convergence of the Lotka-Nagaev estimator.
 In this paper, we establish some Cram\'{e}r moderate deviation results for the Lotka-Nagaev estimator via a martingale method.
 Notice that the Cram\'{e}r moderate deviation results for a supercritical branching process in a random environment (BPRE) stated in  \cite{GLE17} do  not implies our results, because the  random environment for BPRE cannot be degenerate and they considered the estimator $\frac1n \ln Z_n$ instead of the Lotka-Nagaev estimator.

 The paper is organized as follows. In Section \ref{sec2}, we present our main results, including
 Cram\'{e}r   moderate deviations and moderate deviation principles
  for the   Lotka-Nagaev  estimator.
  In Section  \ref{sec3}, we  present some applications of our results in statistics.
 The remaining sections are devoted to the proofs of  theorems.


\section{Main results}\label{sec2}
\setcounter{equation}{0}

\subsection{ Moderate deviations for the weighted  Lotka-Nagaev estimator using the data $(Z_{k})_{n_0\leq k\leq n_0+n}$}
Let $n_0, n \in \mathbb{N}.$ Denote
\begin{eqnarray}\label{fdfmn}
   \widehat{m}_n  = \frac{1}{\sum_{k=n_0}^{n_0+n-1} \sqrt{Z_{ k}} }  \sum_{k=n_0}^{n_0+n-1}    \sqrt{Z_{ k}} (  \frac{Z_{ k+1}}{Z_{ k}} )
\end{eqnarray}
the random weighted  Lotka-Nagaev estimator.
In usual, one takes $n_0=0.$ Here we consider the more general case that $n_0$ may depend on $n$. Denote
$$H_{n_0, n}=\frac{1}{v\sqrt{n}\ } \sum_{k=n_0}^{n_0+n-1}\sqrt{Z_{ k}} \Big( \frac{Z_{ k+1}}{Z_{ k}} -m \Big).$$
 Then $H_{n_0, n}$ can be rewritten in the following form
$$H_{n_0, n}= \frac{ \widehat{m}_n -m   }{v\sqrt{n}   \ }\sum_{k=n_0}^{n_0+n-1}    \sqrt{Z_{ k}} ,$$
and thus it gives a suitable norming of the error term $\widehat{m}_n -m $.
It is easy to check that $$\mathbb{E}\Big[ \sqrt{Z_{ k}} ( \frac{Z_{ k+1}}{Z_{ k}} -m )\Big| Z_0,...,Z_k \Big]=0\ \ \ \textrm{and} \ \ \textrm{Var}\Big(\sqrt{Z_{ k}} ( \frac{Z_{ k+1}}{Z_{ k}} -m )\Big)=v^2.$$
Thus $H_{n_0, n}$ is a standardized martingale, and $(H_{n_0, n})_{n\geq1}$ is the standardized process for the estimator $\widehat{m}_n$.
We have the following Cram\'{e}r moderate deviation result  with respect to $H_{n_0, n}.$
\begin{theorem}\label{th01}
Assume that there exists a positive constant $c$ such that
\begin{eqnarray}\label{Btcondition}
    \mathbb{E}|Z_{ 1}  -m   |^{l } \leq \frac{1}{2} \, l ! \ (l -1 )^{-l/2} \, c^{l-2}\,  \mathbb{E}(Z_1  -m )^{2},\ \ \ \ \ \ l\geq 2.
\end{eqnarray}
Then  the following equalities hold for all $0\leq x = o(  \sqrt{n} ),$
\begin{equation}\label{dfdsf12}
\bigg|\ln\frac{\mathbb{P}(H_{n_0,n} \geq x)}{1-\Phi(x)} \bigg|=O \bigg( \frac{x^3}{\sqrt{n} }+  (1+ x  )\frac{\ln n}{\sqrt{n}}  \bigg)
\end{equation}
and
\begin{equation}
\bigg|\ln\frac{\mathbb{P}(H_{n_0,n} \leq - x)}{ \Phi(-x)} \bigg| =O \bigg( \frac{x^3}{\sqrt{n} }+  (1+ x  )\frac{\ln n}{\sqrt{n}}  \bigg)
\end{equation}
 as $ n \rightarrow \infty.$
\end{theorem}

\begin{remark}
  Let us make some comments on  Theorem \ref{th01}.
\begin{enumerate}


\item It is worth noting that if   $  Z_{ 1} \leq  m+ c_2$, then  condition (\ref{Btcondition}) is satisfied  with
 $c=\frac{1}{3}2^{3/2} \max\{m,\ c_2  \}$.

\item   A sub-Gaussian random variable  also satisfies condition (\ref{Btcondition}), that is,
if there exists a  positive constant $c_1>0$ such that   $$\mathbb{P}(Z_{ 1}  -m \geq x )  \leq  c_1   \exp\{ -x^2/ c_1 \}  , \ \ \  x\geq 0, $$
 then  condition \eref{Btcondition} is satisfied. Indeed, it is easy to see that for all $l\geq 2,$
  \begin{eqnarray}
  \mathbb{E} |Z_{ 1}  -m|^l   &\leq&  m^l \mathbb{P}(Z_{ 1}  -m <0 )  + \int_{0}^\infty l \, x^{l-1} \mathbb{P}(Z_{ 1}  -m \geq x ) dx    \nonumber \\
 &\leq &m^l +\int_{0}^\infty l \, x^{l-1} c_1   \exp\{ -x^2/ c_1 \}  dx= m^l +c_1 (\sqrt{\frac{c_1}{2} }\, )^{l-1} \int_{0}^\infty l \, y^{l-1} \exp\{ -y^2/2 \}  dy \nonumber \\
 &  =& m^l +  c_1 (\sqrt{\frac{c_1}{2} }\, )^{l-1}  l \, !!  .\nonumber
\end{eqnarray}
By Stirling's formula  $n!= \sqrt{2\pi n}\, n^ne^{-n}e^{\frac{1}{12 n\theta_n} } $ for some  $0 \leq \theta_n \leq 1,$
 we deduce that  for all $l \geq 2,$
 \begin{eqnarray}
  l \, !!   &\leq&  \sqrt{(l+1)!\ } \ = \ \frac{l\,! \, (l+1)}{\sqrt{(l+1)!\ } }  \ \leq \ \frac{l\,! \, (l+1)}{\sqrt{ (l+1)^{l+1}e^{-(l+1)}     }\  }
  \ = \ l\,! \, (l+1)^{-(l-1)/2}e^{ (l+1)/2} \nonumber\\
  &\leq&  l\,! \, (l-1)^{-(l-1)/2}e^{ (l+1)/2} \nonumber \\
  &=& l\,! \, (l-1)^{- l /2} (l-1)^{1/2}e^{ (l+1)/2} .\nonumber
\end{eqnarray}
Thus, we have  $$ \mathbb{E} |Z_{ 1}  -m|^l  \leq  m^l +  l\,! \,  (l-1)^{-l/2} c_1 (\sqrt{\frac{c_1}{2} }\, )^{l-1}  (l-1)^{1/2}e^{ (l+1)/2} , $$
which implies (\ref{Btcondition}) with $c$ large enough and $l\geq 3$. When $l=2,$ condition (\ref{Btcondition}) holds obviously.
%

\end{enumerate}
\end{remark}

Using the inequality $|e^x -1 | \leq   e^{ \alpha } |x|$ valid for $|x| \leq \alpha, $ from Theorem  \ref{th01},
we obtain the following result about the equivalence to the normal tail.
\begin{corollary}\label{cor01}
Assume the condition   of Theorem \ref{th01}.
Then for all $n\geq 3$ and all $0\leq x  \leq n^{1/6},$
\begin{equation}\label{dfdsf12}
\frac{\mathbb{P}(H_{n_0,n} \geq x)}{1-\Phi(x)} = 1+ O\bigg( \frac{x^3}{\sqrt{n} }+  (1+ x  )\frac{\ln n}{\sqrt{n}}  \bigg) \ \ \ and \ \ \ \frac{\mathbb{P}(H_{n_0,n} \leq - x)}{ \Phi(-x)} = 1+ O \bigg( \frac{x^3}{\sqrt{n} }+  (1+ x  )\frac{\ln n}{\sqrt{n}}  \bigg).
\end{equation}
In particular, it implies  that
\begin{equation}
 \frac{\mathbb{P}(H_{n_0,n} \geq x)}{1-\Phi(x)} = 1+o(1) \ \ \ and \ \ \ \frac{\mathbb{P}(H_{n_0,n} \leq - x)}{ \Phi(-x)} =1+o(1)
\end{equation}
holds uniformly for  $0 \leq x = o(  n^{1/6}  )$ as $ n \rightarrow \infty.$
\end{corollary}

Theorem  \ref{th01}  also implies  the following  moderate deviation principle  (MDP)  result. For an analogy on a BPRE, but with respect to $\ln Z_n$,
 we refer to Huang and Liu \cite{HT}.
\begin{corollary}\label{co02}
Assume the condition   of Theorem \ref{th01}.
Let $(a_n)_{n\geq1}$ be any sequence of real positive numbers satisfying $a_n \rightarrow \infty$ and $a_n/ \sqrt{n}    \rightarrow 0$
as $n\rightarrow \infty$.  Then  for each Borel set $B$,
\begin{eqnarray}
- \inf_{x \in B^o}\frac{x^2}{2}  \leq   \liminf_{n\rightarrow \infty}\frac{1}{a_n^2}\ln \mathbb{P} \bigg(  \frac{H_{n_0,n}  }{ a_n  }     \in B \bigg)
  \leq \limsup_{n\rightarrow \infty}\frac{1}{a_n^2}\ln \mathbb{P} \bigg(\frac{ H_{n_0,n}  }{ a_n }  \in B \bigg) \leq  - \inf_{x \in \overline{B}}\frac{x^2}{2}   ,   \label{MDP}
\end{eqnarray}
where $B^o$ and $\overline{B}$ denote the interior and the closure of $B$, respectively.
\end{corollary}

The proof of Corollary \ref{co02} is given in Section \ref{sec5}.

%
%
%
%

\subsection{Moderate deviations for the  Lotka-Nagaev estimator using the data $Z_{n}$ and $Z_{n+1}$}

For the Galton-Watson process, it holds
$$\mathbb{E} [ (Z_{n+1}  -m Z_{n})^2 | Z_n] = \mathbb{E} [ ( \sum_{i=1}^{Z_n}  (X_{n,i} -m) )^2 | Z_n] = Z_{n}  \upsilon^2.$$
Thus the following one
$$R_n= \frac{ \sqrt{Z_n} \,}{\upsilon } \Big ( \frac{Z_{n+1}}{Z_{n}} -m \Big)$$
is a normalized process  for the Lotka-Nagaev estimator.
When $Z_1$ satisfies the Cram\'{e}r condition (cf.\ (\ref{fqdasf01})), we have the following Cram\'{e}r moderate deviation result for the normalized Lotka-Nagaev estimator $R_n$.
\begin{theorem}\label{scdsds}
Assume  there exists a constant  $\kappa_0>0$
 such that
\begin{eqnarray}\label{fqdasf01}
\mathbb{E}\exp\{  \kappa_0  Z_1  \}<\infty.
\end{eqnarray}
 Then
\begin{eqnarray}\label{fqdqv01}
\bigg|\ln\frac{\mathbb{P} (R_{n} \geq x)}{1-\Phi(x)} \bigg| = O\Big(  \frac{1+x^{3}}{ \sqrt{n} } \Big)
\end{eqnarray}
holds uniformly for $0\leq x= o(\sqrt{n} ) $ as $n\rightarrow \infty$.
In particular, it implies that
\begin{equation}\label{fqdvf01}
\frac{\mathbb{P} \big( R_n \geq x  \big)}{1-\Phi(x)} =1+ o(1)
\end{equation}
holds uniformly for $0\leq x =o(  n^{1/6} ) $ as $n\rightarrow \infty$.
\end{theorem}

Clearly,  the ranges of validity  for  \eqref{fqdqv01}
 and \eqref{fqdvf01}  coincide with the case of classical Cram\'{e}r moderate deviation result  \cite{Cramer38}.

As $Z_1 \geq 0$, we still have the following Cram\'{e}r moderate deviation result for the normalized Lotka-Nagaev estimator $R_n$ under a weaker moment condition.
\begin{theorem}\label{scddsddss}
Assume  that $ \mathbb{E} Z_1 ^{2+\rho}< \infty$ for some  $  \rho \in (0,  1]$.
 Then
\begin{eqnarray}\label{fghl}
\bigg|\ln\frac{\mathbb{P} (R_{n} \leq - x)}{ \Phi(-x)} \bigg| = O\Big(  \frac{1+x^{2+\rho}}{  n^{\rho/2} } \Big)
\end{eqnarray}
holds uniformly for $0\leq x = o(\sqrt{n}  )$ as $n\rightarrow \infty$.
In particular, it implies that
\begin{equation}
\frac{\mathbb{P} \big( R_n \leq - x  \big)}{ \Phi(-x)} =1+ o(1)
\end{equation}
holds uniformly for $0\leq x=o(  n^{\rho/(4+2\rho) }  )$ as $n\rightarrow \infty$.
\end{theorem}

Clearly, condition (\ref{fqdasf01}) implies that $ \mathbb{E} Z_1 ^{3}< \infty$. Thus, with condition (\ref{fqdasf01}), Theorem \ref{scddsddss} implies that (\ref{fghl})
holds with $\rho=1$. By an argument similar to the
 proof of Corollary \ref{co02}, we have following MDP  result for $R_{n}$.
\begin{corollary}\label{corollary02}
Assume the condition   of Theorem \ref{scdsds}.
Let $(a_n)_{n\geq1}$ be a  sequence of real numbers satisfying $a_n \rightarrow \infty$ and $a_n/ \sqrt{n}    \rightarrow 0$
as $n\rightarrow \infty$.  Then  for each Borel set $B$,
\begin{eqnarray}
- \inf_{x \in B^o}\frac{x^2}{2}  \leq   \liminf_{n\rightarrow \infty}\frac{1}{a_n^2}\ln \mathbb{P} \bigg(  \frac{R_n }{ a_n  }     \in B \bigg)
  \leq \limsup_{n\rightarrow \infty}\frac{1}{a_n^2}\ln \mathbb{P} \bigg(\frac{ R_n }{ a_n }  \in B \bigg) \leq  - \inf_{x \in \overline{B}}\frac{x^2}{2}   ,   \label{MDP}
\end{eqnarray}
where $B^o$ and $\overline{B}$ denote the interior and the closure of $B$, respectively.
\end{corollary}

Under the Linnik condition \cite{L61} (instead of  Cram\'{e}r's condition \eref{fqdasf01}),
we have the following Cram\'{e}r type moderate deviation result  for  the normalized Lotka-Nagaev estimator $R_n$.
\begin{theorem} \label{fsfh25}
Assume that there exist two constants  $\iota_0>0$ and $\tau \in (0, \frac16]$
 such that
\begin{eqnarray}\label{linnik}
\mathbb{E}\exp\{  \iota_0   Z_1^{\frac{4\tau}{ 2\tau+1 }}  \}<\infty .
\end{eqnarray}
Then
\begin{eqnarray}\label{rsfezero}
 \frac{\mathbb{P} (R_{n} \geq x)}{1-\Phi(x)} =1+o(1)
\end{eqnarray}
holds uniformly for  $ x \in  [0, \, o( n^{\tau }  ))$   as $n\rightarrow \infty$.
\end{theorem}

Inequality (\ref{rsfezero}) states that the relative error of normal approximation for $R_{n}$ tends to zero   uniformly for $x \in [0, o(  n^{\tau } )).$
 The range of validity  for (\ref{rsfezero})   coincides  with the Cram\'{e}r moderate deviation of Linnik \cite{L61}   for i.i.d. random variables.

\section{Applications to construction of confidence intervals}\label{sec3}
\setcounter{equation}{0}

\subsection{Case where the data $(Z_{k})_{n_0\leq k\leq n_0+n}$ is observed}
Cram\'{e}r moderate deviation results can be applied to  construction of confidence intervals for $m$. Recall $\widehat{m}_n$ defined by  (\ref{fdfmn}). By Theorem \ref{th01},
we have the following result for the confidence interval  for $m$.
\begin{proposition}\label{c0kldfddgs}
Assume the condition of Theorem \ref{th01}.   Let $\kappa_n \in (0, 1).$  Assume
\begin{eqnarray}\label{keldffdet}
 \big| \ln \kappa_n \big| =o \big( n^{1/3}  \big) .
\end{eqnarray}
Then $[A_{n_0, n},B_{n_0, n}]$, with
\begin{eqnarray*}
A_{n_0, n}=\widehat{m}_n - \frac{  v \sqrt{n}\, \Phi^{-1}( 1-\kappa_n/2) }{ \sum_{k=n_0}^{n_0+n-1}\sqrt{Z_{ k}}  }
 \ \ \
\textrm{and}
\ \ \
B_{n_0, n}=\widehat{m}_n + \frac{  v \sqrt{n}\, \Phi^{-1}( 1-\kappa_n/2) }{ \sum_{k=n_0}^{n_0+n-1}\sqrt{Z_{ k}}  },
\end{eqnarray*}
is a  $1-\kappa_n$ confidence interval for $m$, for $n$ large enough.
\end{proposition}
\emph{Proof.}  Notice that $ 1-\Phi \left( x\right) =  \Phi \left( -x\right). $ Corollary \ref{cor01} implies that
\begin{equation}   \label{sfdfsdfs}
\frac{\mathbb{P}( H_{n_0,n}  \geq x)}{1-\Phi \left( x\right)}=1+o(1)\ \ \ \ \textrm{and} \ \ \ \ \frac{\mathbb{P}(H_{n_0,n}  \leq- x)}{ \Phi \left(- x\right)}=1+o(1)
\end{equation}
uniformly for $0\leq x=o (   n^ {1/6}   )$.
Notice that the  inverse function $\Phi^{-1}$
of the the standard normal distribution function $\Phi$ satisfies the following asymptotic expansion
$$\Phi^{-1}(1-p_n)=\sqrt{\ln(1/p_n^2)-\ln\ln(1/p_n^2) -\ln(2\pi) }+ o(p_n) ,\ \ \ \ \ p_n \searrow 0. $$
By (\ref{fgsgjds1}) and \eref{keldffdet}, it is easy to see that the upper $(\kappa_n/2)$th quantile of a standard normal distribution
$\Phi^{-1}( 1-\kappa_n/2)=-\Phi^{-1}(   \kappa_n/2) = O(\sqrt{| \ln \kappa_n |} ) $   is of order $o\big( n^{1/6}  \big).$
Then applying the last equality to (\ref{sfdfsdfs}), we have
\[
\mathbb{P}\big(H_{n_0,n}   \geq  \Phi^{-1}( 1-\kappa_n/2)\big) \sim \kappa_n/2
\ \ \ \
\textrm{and}
\ \ \ \
  \mathbb{P}\big( H_{n_0,n}  \leq -\Phi^{-1}( 1-\kappa_n/2) \big) \sim \kappa_n/2
\]
as $n\rightarrow \infty.$
Clearly, $H_{n_0,n}  \leq \Phi^{-1}( 1-\kappa_n/2)$ means that $m \geq A_{n_0, n},   $ while $H_{n_0,n}  \geq -\Phi^{-1}( 1-\kappa_n/2)$
means  $m\leq B_{n_0, n}.  $ This completes the proof of  Proposition \ref{c0kldfddgs}. \hfill\qed

\subsection{Case where the data $Z_{n}$ and $Z_{n+1}$ is observed}
When $Z_{n+1}$ and $Z_{n}$ can be observed, we can make use of  Theorem  \ref{fsfh25}
to construct confidence intervals.
\begin{proposition}\label{c0kls}
Assume the condition of Theorem \ref{fsfh25}. Let $\kappa_n \in (0, 1).$  Assume
\begin{eqnarray}\label{keldetd}
 \big| \ln \kappa_n \big| =o \big( n^{2\tau }  \big) .
\end{eqnarray}
Let $$\Delta_n=\frac{\Phi^{-1}(1-\kappa_n/2)  }{Z_n } v.$$
Then   $[A_n,B_n]$, with
\begin{eqnarray*}
A_n=\frac{ Z_{n+1} }{ Z_{n}}  -\Delta_n  \quad   \textrm{and} \ \quad  
B_n=\frac{ Z_{n+1} }{ Z_{n}} +\Delta_n, 
\end{eqnarray*}
is a  $1-\kappa_n$ confidence interval for $m $, for $n$ large enough.
\end{proposition}
\emph{Proof.} Theorem  \ref{fsfh25} implies that
\begin{equation} \label{tphisnds4}
\frac{\mathbb{P} (R_{n}  \geq x)}{1-\Phi \left( x\right)}=1+o(1)\ \ \ \ \textrm{and} \ \ \ \ \frac{\mathbb{P} (R_{n}  \leq- x)}{1-\Phi \left( x\right)}=1+o(1)
\end{equation}
uniformly for $0\leq x=o (   n^{\tau }  ).$  When $\kappa_n$ satisfies the condition (\ref{keldetd}),
  the upper $(\kappa_n/2)$th quantile of a standard normal distribution satisfies
 $\Phi^{-1}( 1-\kappa_n/2) = O(\sqrt{| \ln \kappa_n |} ), $  which  is of order $o\big( n^\tau  \big).$
Using (\ref{tphisnds4}), by an argument similar to the
 proof of  Proposition \ref{c0kldfddgs},  we obtain the $1-\kappa_n$ confidence interval for $m.$  \hfill\qed

 \section{Proof of Theorem \ref{th01}}
 \setcounter{equation}{0}
  Let $(\xi_i, \mathcal{F}_i)_{1\leq i \leq n}$
be a finite sequence of martingale differences.
In the sequel we shall use the following conditions:
\begin{description}
\item[(A1)]  There exists a number $\epsilon_n \in (0, \frac12]$ such that
\[
|\mathbb{E}[\xi_{i}^{k}  | \mathcal{F}_{i-1}]| \leq \frac12 k!\epsilon_n^{k-2} \mathbb{E}[\xi_{i}^2 | \mathcal{F}_{i-1}],\ \ \ \ \ \textrm{for all}\ k\geq 3\ \ \textrm{and}\ \ 1\leq i\leq n;
\]
\item[(A2)]  There exists a number  $ \delta_n \in [0, \frac12]$ such that
$ \left| \sum_{i=1}^n\mathbb{E}[\xi_{i}^2 | \mathcal{F}_{i-1}]-1\right| \leq  \delta_n^2.$
\end{description}

 In the proof of Theorem \ref{th01}, we make use of the following lemma  which gives a Cram\'{e}r moderate deviation result for martingales. See Theorems 2.1 and 2.2 of Fan, Grama and Liu \cite{FGL13}.
\begin{lemma}\label{th0}
Assume conditions (A1) and (A2). Then there exits an absolute constant $\alpha \in (0,1)$ such that for all $0\leq x \leq \alpha \, \epsilon_n^{-1}$
and $\delta_n \leq \alpha$,
\begin{equation}\nonumber
\Bigg| \ln \frac{\mathbb{P}(\sum_{i=1}^n \xi_i \geq x)}{1-\Phi \left( x\right)} \Bigg| \leq   C_{\alpha}  \bigg(x^3 \epsilon_n + x^2 \delta_n^2 +  \,  (1+ x)\left( \epsilon_n \left| \ln  \epsilon_n
 \right| + \delta_n \right) \bigg)
\end{equation}
and
\begin{equation} \nonumber
\Bigg| \ln \frac{\mathbb{P}(\sum_{i=1}^n \xi_i \leq - x)}{ \Phi \left(- x\right)} \Bigg| \leq   C_{\alpha}  \bigg(x^3 \epsilon_n + x^2 \delta_n^2 +  \,  (1+ x)\left( \epsilon_n \left| \ln  \epsilon_n
 \right| + \delta_n \right) \bigg) ,
\end{equation}
where the constant $C_{\alpha}$ does not depend on $(\xi _i,\mathcal{F}_i)_{i=0,...,n}$, $n$ and $x$.
\end{lemma}

Denote $$\hat{\xi}_{k+1}= \sqrt{Z_{ k}} ( Z_{ k+1}/Z_{ k} -m ),$$
$\mathfrak{F}_{n_0} =\{ \emptyset, \Omega \}  $ and $\mathfrak{F}_{k+1}=\sigma \{ Z_{i}: n_0\leq i\leq k+1  \}$ for all $k> n_0$.
Notice that $X_{k,i}$ is independent of $Z_k.$
Then it is easy to verify that $\mathbb{E}[ \hat{\xi}_{k+1}  |\mathfrak{F}_{k } ] =0.$
 Thus
 $(\hat{\xi}_k, \mathfrak{F}_k)_{k=n_0+1,...,n_0+n}$ is a finite sequence of martingale  differences.
Notice that   $X_{ k, i}-m, i\geq 1,$  are  centered and independent random variables. Thus,
 the following   equalities hold
\begin{eqnarray}
  \mathbb{E}[ \hat{\xi}_{k+1}^2  |\mathfrak{F}_{k } ]
 =    v^2  \ \ \ \textrm{and} \ \ \ \sum_{k=n_0}^{n_0+n-1} \mathbb{E}[ \hat{\xi}_{k+1}^2  |\mathfrak{F}_{k } ]
 =  n v^2 .  \label{ineq9.2}
\end{eqnarray}
By Rio's inequality (cf. Theorem 2.1 of Rio \cite{R09}) and the fact that $X_{ k, i}$ is independent to $\mathfrak{F}_{k }$, we have
for any $  \rho \geq 1,$
\begin{eqnarray*}
\Big (\mathbb{E}[  | \sum_{i=1}^{Z_k}    (X_{ k, i}  -m)  |^{2\rho }  |\mathfrak{F}_{k } ] \Big)^{2/(2 \rho)}
&\leq& (2\rho -1 )    \sum_{i=1}^{Z_k} \Big(\mathbb{E}| X_{ k, i}  -m   |^{2\rho } \Big)^{2/(2\rho)} .
\end{eqnarray*}
The last inequality implies that for any $  \rho \geq 1,$
\begin{eqnarray*}
\mathbb{E}[  | \sum_{i=1}^{Z_k}    (X_{ k, i}  -m)  |^{2\rho }  |\mathfrak{F}_{k } ]
&\leq& (2\rho -1 )^\rho   \Big( \sum_{i=1}^{Z_k} \Big(\mathbb{E}| X_{ k, i}  -m   |^{2\rho } \Big)^{1/ \rho}  \Big)^\rho \\
&\leq& (2\rho -1 )^\rho Z_k^\rho     \mathbb{E}| X_{ k, 1}  -m   |^{2\rho }  .
\end{eqnarray*}
Hence, the following inequalities hold for any $  \rho \geq 1,$
\begin{eqnarray}
  \mathbb{E}[ |\hat{\xi}_{k+1}|^{2 \rho} |\mathfrak{F}_{k } ] &=&
  Z_{k}^{- \rho } \mathbb{E}[  | Z_{ k+1} -mZ_{ k} |^{2 \rho }  |\mathfrak{F}_{k } ]
 =    Z_{k}^{ -\rho }   \mathbb{E}[  | \sum_{i=1}^{Z_k}    (X_{ k, i}  -m)  |^{2 \rho }  |\mathfrak{F}_{k } ] \nonumber\\
 & \leq & (2\rho -1 )^\rho  \mathbb{E}| X_{ k, 1}  -m   |^{2\rho } . \nonumber
\end{eqnarray}
The last inequality becomes equality when $\rho=1$.
Notice that $X_{ k, 1}$ has the same distribution as $Z_1.$
Thus, by condition (\ref{Btcondition}), we get for all $l\geq 2,$
\begin{eqnarray}
  \mathbb{E}[ |\hat{\xi}_{k+1}|^{l} |\mathfrak{F}_{k } ] &\leq &    (l -1 )^{l/2}     \mathbb{E}|X_{ k, 1}  -m   |^{l }
   \leq (l -1 )^{l/2}  \frac{1}{2} \, l ! \ (l -1 )^{-l/2} \, c^{l-2}\,  \mathbb{E}(X_{ k, 1}  -m )^{2} \nonumber \\
  &=&   \frac12   \, l ! \, c^{l-2}\,  v^2
    \ =\    \frac12   \, l ! \, c^{l-2}\,  \mathbb{E}[ \hat{\xi}_{k+1}^2  |\mathfrak{F}_{k } ] .    \nonumber
\end{eqnarray}
Set $\xi_k= \hat{\xi}_{n_0+k }/\sqrt{n} v$ and $\mathcal{F}_{k}=\mathfrak{F}_{n_0+k }$. It is easy to see that conditions (A1) and (A2) are satisfied with
$\epsilon_n=c/ \sqrt{n} v$ and $\delta_n=0$. Applying  Lemma \ref{th0} to $(\xi_k, \mathcal{F}_k)_{1\leq k \leq n}$, we obtain the desired inequalities.

\section{Proof of Corollary \ref{co02}}\label{sec5}
\setcounter{equation}{0}
We first show that for any Borel set $B\subset \mathbb{R},$ it holds
\begin{eqnarray}\label{dfgsfdd}
\limsup_{n\rightarrow \infty}\frac{1}{a_n^2}\ln \mathbb{P}\bigg(\frac{H_{n_0,n} \ }{ a_n  }  \in B \bigg)   \leq  - \inf_{x \in \overline{B}}\frac{x^2}{2}.
\end{eqnarray}
When $B  =\emptyset,$ the last inequality holds with   $-\inf_{x \in \emptyset}\frac{x^2}{2}=-\infty$. Hence,  we only need to consider the case of $B  \neq \emptyset.$ Let $x_0=\inf_{x\in B} |x|,$ then we have $x_0=\inf_{x\in \overline{B}} |x|.$
Then, from  Theorem \ref{th01}, it follows that for $a_n =o(\sqrt{n}),$
\begin{eqnarray*}
 \mathbb{P}\bigg(\frac{ H_{n_0,n} \ }{ a_n  }  \in B \bigg)
 &\leq&  \mathbb{P}\bigg( | H_{n_0,n}|  \geq a_n x_0\bigg)\\
 &\leq&  2\Big( 1-\Phi \left( a_nx_0\right)\Big)
  \exp\bigg\{ C   \bigg(  \frac{( a_nx_0)^3}{\sqrt{n} }+  (1+ ( a_nx_0)  )\frac{\ln n}{\sqrt{n}}       \bigg)  \bigg\}.
\end{eqnarray*}
Using the following two-sided bound for the normal distribution function
\begin{eqnarray}\label{fgsgjds1}
\frac{1}{\sqrt{2 \pi}(1+x)} e^{-x^2/2} \leq 1-\Phi ( x ) \leq \frac{1}{\sqrt{ \pi}(1+x)} e^{-x^2/2}, \ \   x\geq 0,
\end{eqnarray}
  and  the fact that $a_n \rightarrow \infty$ and $a_n/\sqrt{n}\rightarrow 0$,
we obtain
\begin{eqnarray*}
\limsup_{n\rightarrow \infty}\frac{1}{a_n^2}\ln \mathbb{P}\bigg(\frac{ H_{n_0,n}   \ }{ a_n  }  \in B \bigg)
 \ \leq \  -\frac{x_0^2}{2} \ = \  - \inf_{x \in \overline{B}}\frac{x^2}{2} ,
\end{eqnarray*}
which gives (\ref{dfgsfdd}).

Next, we show that the following inequality holds
\begin{eqnarray}\label{dfsfffgsfn}
\liminf_{n\rightarrow \infty}\frac{1}{a_n^2}\ln \mathbb{P}\bigg(\frac{ H_{n_0,n} \ }{ a_n   }  \in B \bigg) \geq   - \inf_{x \in B^o}\frac{x^2}{2} .
\end{eqnarray}
When $B^o =\emptyset,$ the last inequality holds obvious, with   $ -\inf_{x \in  \emptyset}\frac{x^2}{2}=-\infty$. Hence, we may assume that $B^o \neq \emptyset.$ Since $B^o$ is an open set,
for any given small positive constant $\varepsilon_1,$ there exists an $x_0 \in B^o,$ such that
\begin{eqnarray*}
 0< \frac{x_0^2}{2} \leq   \inf_{x \in B^o}\frac{x^2}{2} +\varepsilon_1.
\end{eqnarray*}
By the fact that $B^o$ is an open set again, for $x_0 \in B^o$ and any  $\varepsilon_2 \in (0, |x_0|], $ it holds $(x_0-\varepsilon_2, x_0+\varepsilon_2]  \subset B^o.$  Without loss of generality, we may assume that $x_0>0.$
Then, we have
\begin{eqnarray}
\mathbb{P}\bigg(\frac{H_{n_0,n} \ }{ a_n  } \in B  \bigg)   &\geq&   \mathbb{P}\bigg( H_{n_0,n}    \in (a_n ( x_0-\varepsilon_2), a_n( x_0+\varepsilon_2)] \bigg) \nonumber \\
&=&   \mathbb{P}\bigg( H_{n_0,n} \geq  a_n ( x_0-\varepsilon_2)   \bigg)-\mathbb{P}\bigg( H_{n_0,n}  \geq  a_n( x_0+\varepsilon_2) \bigg). \label{fsfvhd}
\end{eqnarray}
By Theorem  \ref{th01}, it is easy to see that for $a_n \rightarrow \infty$ and $ a_n =o(\sqrt{n} ),$ $$\lim_{n\rightarrow \infty} \frac{\mathbb{P}\big(H_{n_0,n} \geq  a_n( x_0+\varepsilon_2) \big) }{\mathbb{P}\big( H_{n_0,n} \geq  a_n ( x_0-\varepsilon_2)   \big)  } =0 .$$
By the last line and Theorem \ref{th01}, it holds for all $n$ large enough and $a_n =o(\sqrt{n} ),$
\begin{eqnarray*}
&&\mathbb{P}\bigg(\frac{H_{n_0,n} \ }{ a_n  } \in B  \bigg) \ \geq\    \frac12 \mathbb{P}\bigg(  H_{n_0,n}   \geq  a_n ( x_0-\varepsilon_2)   \bigg) \\
&& \ \ \ \ \ \ \ \ \  \ \geq   \   \frac12 \Big( 1-\Phi \left( a_n( x_0-\varepsilon_2)\right)\Big)   \exp\bigg\{ -C  \bigg(  \frac{( a_n( x_0-\varepsilon_2))^3}{\sqrt{n} }+  (1+   a_n( x_0-\varepsilon_2)   )\frac{\ln n}{\sqrt{n}}       \bigg) \bigg\}.
\end{eqnarray*}
Using   (\ref{fgsgjds1}) and  the fact that $a_n \rightarrow \infty$ and $a_n/\sqrt{n}\rightarrow 0$, after some simple calculations,
we obtain
\begin{eqnarray*}
 \liminf_{n\rightarrow \infty}\frac{1}{a_n^2}\ln \mathbb{P}\bigg(\frac{H_{n_0,n}   \ }{ a_n  }  \in B \bigg)  \geq  -  \frac{1}{2}( x_0-\varepsilon_2)^2 . \label{ffhms}
\end{eqnarray*}
Letting $\varepsilon_2\rightarrow 0,$  we  arrive at
\begin{eqnarray*}
\liminf_{n\rightarrow \infty}\frac{1}{a_n^2}\ln \mathbb{P}\bigg(\frac{ H_{n_0,n}  \ }{ a_n  }  \in B \bigg) \ \geq\ -  \frac{x_0^2}{2}  \  \geq \   -\inf_{x \in B^o}\frac{x^2}{2} -\varepsilon_1.
\end{eqnarray*}
Since that $\varepsilon_1$ can be arbitrarily small, we get (\ref{dfsfffgsfn}).
Combining (\ref{dfgsfdd}) and (\ref{dfsfffgsfn}) together, we obtain the desired result. This  completes the proof of Corollary  \ref{co02}.
 \hfill\qed

\section{Proof  of Theorem  \ref{scdsds} }
\setcounter{equation}{0}
\subsection{Preliminary lemmas}
Denote by $f_{n}(s)=\mathbb{E} s^{Z_{n}},$ $ |s|\leq 1,$ the generating function of $Z_{n}$.
In the proof  of Theorem  \ref{scdsds}, we shall make use of  the following lemma, see Athreya \cite{A94}.
\begin{lemma}\label{lemma01}
If  $ p_1>0$, then it holds
\begin{equation}
\lim\limits_{n\to\infty}\frac{f_n(s)}{p_1^n}=\sum_{k=1}^{\infty}q_{k}s^k,
\end{equation}
where $(q_{k}, k\geq 1)$   is defined by the generating function $Q(s)=\sum_{k=1}^{\infty}q_{k}s^k, 0\leq s<1,$ the unique solution of the following functional equation
$$Q(f(s))=p_{1}Q(s), \quad \mbox{ where\ \  } f(s) =\sum_{j=1}^{\infty}p_{j}s^j, \ 0\leq s<1,$$
subject to
$$Q(0)=0, \qquad  Q(1)=\infty, \qquad Q(s)<\infty \mbox{\ \ for\ \ } 0\leq s <1.$$
\end{lemma}

By Lemma \ref{lemma01}, we obtain the following estimation for $Z_n$.
\begin{lemma}\label{lemma2}
It holds
\begin{eqnarray}\label{fgffdgdfg}
\mathbb{P}(Z_n  \leq n) \leq   C_1 \exp\{  -n c_0\}.
\end{eqnarray}
\end{lemma}
\textit{Proof.}
When $p_1>0,$ by Markov's inequality,  we deduce that for $s_0=\frac{1+p_1}{2} \in (0, 1),$
\begin{eqnarray*}
 \mathbb{P}(Z_n  \leq n) = \mathbb{P}(s_0^{Z_n}\geq s_0^n) \leq  s_0^{-n}f_n(s_0).
\end{eqnarray*}
Using Lemma \ref{lemma01},  we have
\begin{eqnarray}
 \mathbb{P}(Z_n  \leq n)
&\leq& C (\frac{p_1}{s_0})^nQ(s_0) \nonumber \\
&=& C_1 \exp\{  -n c_0\},
\end{eqnarray}
where $C_1=C  Q(s_0) $ and $c_0= \ln (s_0/p_1).$ Notice that
$s_0  \in (p_1, 1),$ which implies that $c_0>0.$
When $p_1=p_0=0,$   we have $Z_n\geq 2^n,$ and (\ref{fgffdgdfg}) holds obviously for all $n$.
 \hfill\qed

In the proof of Theorem \ref{scdsds}, we
also make use of the  following lemma of Cram\'{e}r \cite{Cramer38}.
\begin{lemma}\label{Cramerlema}
Let $(X_{i})_{i\geq1} $ be i.i.d.\ and centered random variables. Assume that $\mathbb{E}\exp\{  \lambda   |X_1| \}<\infty $ for some constant $ \lambda >0.$ Set $S_{n}=\sum_{i=1}^{n} X_{i}$ and $v^{2}=\mathbb{E}X_{1}^2$.
Then
\begin{equation}
\bigg|\ln\frac{\mathbb{P}(S_{n}/  (v \sqrt{n}) \geq x)}{1-\Phi(x)} \bigg| \leq  C  \frac{1+x^{3}}{ \sqrt{n} }
\end{equation}
uniformly for $0 \leq x =o(\sqrt{n} )$.
\end{lemma}

\subsection{Proof of Theorem  \ref{scdsds} }
By  the definition of $R_n$,  it is easy to see that $R_n$ can be rewritten as follows:
$$R_n =   \frac{ 1}{\upsilon \sqrt{Z_n}} \Big (  Z_{n+1}  -m  Z_{n}\Big)= \frac{ 1}{\upsilon \sqrt{Z_n}}  \sum_{i=1}^{Z_n} \Big ( X_{n,i}  - m  \Big).$$
By the total probability formula, we have
\begin{eqnarray}
\mathbb{P} \Big(R_n \geq  x\Big)
&=&\sum_{k=1}^{\infty}\mathbb{P} (Z_n=k)\mathbb{P} \bigg( \frac{ 1}{\sqrt{k} \upsilon} \sum_{i=1}^{k} \Big ( X_{n,i}  - m  \Big)  \geq x  \bigg)\nonumber \\
&=:&\sum_{k=1}^{\infty}\mathbb{P}(Z_n=k) F_{k}(x).   \label{sfsfds01}
\end{eqnarray}
Notice that $X_{n,i}, \ 1\leq i\leq k,$ have the same distribution as $Z_1$ and $Z_1\geq 0$.
By (\ref{fqdasf01}), it holds
$$ \mathbb{E}e^{  \kappa_0 | Z_1-m|  }  \leq   e^{\kappa_0  m }+  e^{- \kappa_0  m }\mathbb{E}e^{  \kappa_0  Z_1  }<\infty. $$
When $k\geq n,$ by condition \eqref{fqdasf01} and  Lemma \ref{Cramerlema}, we get
\begin{equation}\label{fsfdfs55}
\bigg|\ln\frac{F_{k}(x)}{1-\Phi(x)} \bigg| \leq  C_1  \frac{1+x^{3}}{ \sqrt{k} } \leq  C_1  \frac{1+x^{3}}{ \sqrt{n} }
\end{equation}
uniformly for $0 \leq x =o(\sqrt{n} )$.
Returning to (\ref{sfsfds01}), by (\ref{fsfdfs55}), we have for all $0 \leq x =o(\sqrt{n} )$,
\begin{eqnarray}
\mathbb{P} \Big(R_n \geq   x\Big) &\geq&  \sum_{k=n}^{\infty}\mathbb{P}(Z_n=k) F_{k}(x)
 \geq \Big(  1-\Phi(x) \Big)  \exp\Big\{ - C_1  \frac{1+x^{3}}{ \sqrt{n} } \Big\}    \sum_{k=n}^{\infty}\mathbb{P}(Z_n=k) \nonumber \\
&\geq& \Big(  1-\Phi(x) \Big)  \exp\Big\{ - C_1  \frac{1+x^{3}}{ \sqrt{n} } \Big\}  \Big( 1-  \mathbb{P}(Z_n \leq n)\Big) .  \label{ineq94fg}
\end{eqnarray}
By Lemma \ref{lemma2}, we have
\begin{eqnarray}\label{fssvgsf}
 \mathbb{P}(Z_n\leq n)   \leq C_2 \exp\{  -C_3 n  \}.
 \end{eqnarray}
Applying the last inequality  to \eqref{ineq94fg}, we obtain for all $0 \leq x =o(\sqrt{n} )$,
\begin{eqnarray}\label{fsvrfv}
\mathbb{P} \Big(R_n \geq x\Big)  &\geq& \Big(  1-\Phi(x) \Big)  \exp\Big\{ - C_1  \frac{1+x^{3}}{ \sqrt{n} } \Big\} \Big( 1-  C_2\exp\{  -C_3 n  \}\Big) \nonumber \\
 &\geq&  \Big(  1-\Phi(x) \Big)  \exp\Big\{ - C_4  \frac{1+x^{3}}{ \sqrt{n} } \Big\}.
\end{eqnarray}
Returning to (\ref{sfsfds01}), by (\ref{fsfdfs55}) and \eqref{fssvgsf}, we  deduce that  for all $0 \leq x =o(\sqrt{n} )$,
\begin{eqnarray}
\mathbb{P} \Big(R_n\geq  x\Big) &\leq& \sum_{k=1}^{n-1}\mathbb{P}(Z_n=k)F_{k}(x)   + \sum_{k=n}^{\infty}\mathbb{P}(Z_n=k)F_{k}(x)\nonumber\\
&\leq&\mathbb{P}(Z_n \leq n-1 )   + \Big(  1-\Phi(x) \Big)  \exp\Big\{   C_1  \frac{1+x^{3}}{ \sqrt{n} } \Big\} \sum_{k=n}^{\infty}\mathbb{P}(Z_n=k)\nonumber\\
&\leq& C_2 \exp\{  -C_3 n  \}  +  \Big(  1-\Phi(x) \Big)  \exp\Big\{   C_1  \frac{1+x^{3}}{ \sqrt{n} } \Big\} \nonumber \\
&\leq&    \Big(  1-\Phi(x) \Big)  \exp\Big\{   C_4  \frac{1+x^{3}}{ \sqrt{n} } \Big\}. \label{fsvdsdffv}
\end{eqnarray}
Combining (\ref{fsvrfv}) and (\ref{fsvdsdffv}) together, we obtain the desired inequality, that is (\ref{fqdqv01}).

\section{Proofs  of Theorems  \ref{scddsddss} and \ref{fsfh25}}
\setcounter{equation}{0}
The proof of Theorem \ref{scddsddss} is similar to   the proof of Theorem \ref{scdsds}. However, instead of using Lemma \ref{Cramerlema},
we should make use of the  following lemma of Fan \cite{F17}.
\begin{lemma}
Let $(X_{i})_{i\geq1} $ be i.i.d.\ and centered random variables. Assume that $ X_1 \geq -C  $ and $\mathbb{E}|X_1|^{2+\rho } <\infty $ for some constants $C>0$ and $ \rho \in (0, 1].$ Let $S_{n}=\sum_{i=1}^{n} X_{i}$ and $\upsilon^{2}=\mathbb{E}X_{1}^2$.
Then
\begin{equation}\nonumber
\bigg|\ln\frac{\mathbb{P}(S_{n}/  (\upsilon \sqrt{n}) \leq - x)}{ \Phi(-x)} \bigg| \leq  C  \frac{1+x^{2+\rho}}{  n^{ \rho/2} }
\end{equation}
holds uniformly for $0 \leq x =o(\sqrt{n} )$.
\end{lemma}

The proof of Theorem \ref{fsfh25} is analogous  to the proof of Theorem \ref{scdsds}. However, instead of using Lemma \ref{Cramerlema},
we should make use of the  following lemma of Linnik \cite{L61}.
\begin{lemma}
Let $(X_{i})_{i\geq1} $ be i.i.d.\ and centered random variables. Assume that $\mathbb{E}\exp\{  \iota_0|X_1|^{\frac{4\tau }{2\tau+1}} \}<\infty $
for two constants $\iota_0>0$ and $\tau \in (0, \frac16]$.
 Let $S_{n}=\sum_{i=1}^{n} X_{i}$ and $\upsilon^{2}=\mathbb{E}X_{1}^2$.
Then
\begin{equation}
 \frac{\mathbb{P}(S_{n}/  (\upsilon \sqrt{n}) \geq x)}{1-\Phi(x)}=1+o(1) \nonumber
\end{equation}
holds uniformly for $0 \leq x =o(n^{\tau})$.
\end{lemma}

%

\subsection*{Acknowledgements}
We would like to thank the two referees for their helpful remarks and suggestions.
This work has been partially supported by the National Natural Science Foundation
of China (Grant No. 11971063), by CY
(Grant "Investissements d'Avenir" ANR-16-IDEX-0008)  "EcoDep" PSI-AAP2020-0000000013) and  Labex MME-DII
 (ANR-11-LABEX-0023-01).

\end{document}